\documentclass[12pt,a4paper]{article}
\usepackage{amsmath,amssymb,amsthm,amsfonts,amscd,euscript,verbatim}
\usepackage{t1enc}
\usepackage[cp1251]{inputenc}

\theoremstyle{plain}
\newtheorem{theo}{Theorem}[subsection]
\newtheorem*{theor}{Theorem}
\newtheorem{pr}[theo]{Proposition}

 \newtheorem{lem}[theo]{Lemma}

\theoremstyle{remark}
\newtheorem{rema}[theo]{Remark}

\theoremstyle{remark}
\newtheorem{defi}[theo]{Definition}

\newtheorem*{notata}{Notation and conventions}

 \newcommand\de{\Delta}

\newcommand\al{\alpha}

\newcommand\ok{{\mathfrak{O}_K}}

\newcommand\ol{{\mathfrak{O}_L}}

\newcommand\ovl{{{\overline L}}}
\newcommand\olv{{{\overline l}}}
\newcommand\mw{\tilde{M}{}}

\newcommand\gm{\mathfrak{M}}
\newcommand\ga{\mathfrak{A}}

\newcommand\z{\mathbb{Z}}

\newcommand\la{\lambda}

\newcommand\ob{^{-1}} \newcommand\lan{\langle}
\newcommand\ra{\rangle}
\newcommand\nde{\triangleleft}
\newcommand\pdp{\pi D_{F_\pi}}
\newcommand\ns{\{0\}}
\newcommand\bff{\mathbf{f}}
\newcommand\bv{\mathbf{V}}

\DeclareMathOperator\homm{\operatorname{Hom}}
\DeclareMathOperator\ext{\operatorname{Ext}}
\DeclareMathOperator\id{\operatorname{id}}
\DeclareMathOperator\ke{\operatorname{Ker}}
\DeclareMathOperator\imm{\operatorname{Im}}
\DeclareMathOperator\cok{\operatorname{Coker}}

\DeclareMathOperator\car{\operatorname{Cart}}

\DeclareMathOperator\cl{\operatorname{Cl}}
\DeclareMathOperator\spe{\operatorname{Spec}}

\DeclareMathOperator\prli{\varprojlim}

\begin{document}
 \title{Finite flat commutative group schemes over complete discrete valuation
 rings II: classification, tangent spaces, and semistable reduction of Abelian varieties}
 \author{M.V. Bondarko
\thanks{Supported by Russian Fundamental Science Foundation, grant No. 04-01-00082a. The author is deeply grateful to prof. Ju. Tschinkel and  the
Mathematical Department of the University of G\"ottingen for
providing excellent working conditions. }}
 \date{ December 2004}
\begin{abstract}

Results of  previous papers are used to obtain a complete
classification of finite local flat commutative group schemes over
mixed characteristic complete discrete valuation rings. We
classify group schemes in terms of their Cartier modules. We also
prove the equivalence of different definitions of the tangent
space and the dimension for these group schemes. In particular we
prove that the minimal dimension of a formal group law that
contains a given local group scheme $S$ as a closed subgroup is
equal to the minimal number of
  generators for the coordinate ring of $S$.
As an application the following reduction criteria
 for Abelian varieties are proved.

Let $K$ be a mixed characteristic local field, let its residue
field have characteristic $p$,
 $L$ be a finite  extension of
$K$, let $\mathfrak{O}_K\subset\mathfrak{O}_L$ be their rings of
integers. Let $e$ be the absolute ramification index of $L$,
$s=[\log_p(pe/(p-1))]$, $e_0$ be the ramification index of $L/K$,
$l=2s+v_p(e_0)+1$.

For a finite flat commutative $\mathfrak{O}_L$-group scheme $H$ we
denote the $\mathfrak{O}_L$-dual of the module
 $J/J^2$ by $TH$. Here $J$ is the augmentation ideal of the
coordinate ring of $H$.

Let $V$ be an $m$-dimensional Abelian variety over $K$. Suppose
that $V$ has semistable reduction over $L$.

\begin{theor}[A]

$V$ has semistable reduction over $K$ if and only if
 for some group scheme $H$ over $\mathfrak{O}_K$
 there exist embeddings of $H_K$ into
$\operatorname{Ker}[p^{l}]_{V,K}$, and of
$(\mathfrak{O}_L/p^l\mathfrak{O}_L)^m$ into $TH_\ol$.
\end{theor}

This criterion has a very nice-looking version in the ordinary
reduction case.

\begin{theor}[B]
   $V$ has ordinary  reduction over $K$
if and only if for some  $H_K\subset
\operatorname{Ker}[p^{l}]_{V,K}$ and $M$ unramified
  over $K$ we have
  $H_M\cong (\mu_{p^{l},M})^m$.
Here $\mu$ denotes the group scheme of roots of unity.
\end{theor}


\end{abstract}


 \maketitle

\section*{Introduction}

In the paper \cite{02}  the Cartier modules of finite flat local
commutative  group schemes over the rings of integers of complete
discrete valuation fields (defined by Oort) were used to obtain
results on the generic fibre of (finite flat commutative) group
schemes. As an application a  certain finite (see below) $p$-adic
good reduction criterion for Abelian varieties was proved.

In this paper we use  results  of the previous work and obtain a
complete classification of finite local flat commutative group
schemes over mixed characteristic complete discrete valuation
rings in terms of their Cartier modules. We give an explicit
description of the image of the Oort functor. The classification
obtained has the following advantages when compared with the
classification of Breuil (see \cite{clb}).

1. The modules that appear as $C(S)$ are explicitly described for
arbitrary local group schemes. In \cite{clb} an explicit
classification is given only for $S$ satisfying certain extra
flatness conditions; the modules for arbitrary $S$ are obtained by
extensions. We note that the subclass of group schemes that was
classified in \cite{clb} explicitly  also can be easily described
(i.e. distinguished from other local group schemes) in  terms of
their Cartier modules.

2. The Cartier module language is Galois-stable; this is crucial
for descent questions.

3. We don't restrict ourselves to the perfect residue field case.

On the other hand, the advantage of \cite{clb} is that  $p$-group
schemes there were not assumed to be local.

Our classification will probably be generalized to not necessarily
local $p$-group schemes in a succeeding paper.

We also prove the equivalence of different definitions of the
tangent spaces and the dimension of these group schemes. As an
application  certain finite $p$-adic criteria for semistable and
ordinary reduction of Abelian varieties are proved. We call these
criteria  finite because in contrast to Grothendieck's criteria
(see \cite{gro}) it is sufficient to check certain conditions on
some finite $p$-torsion subgroup of $V$ (instead of the whole
$p$-torsion).

We introduce some notation.

Let $K$ be a complete discrete valuation field of characteristic
$0$ with residue field of characteristic
 $p$, let
 $L$ be a finite extension of
$K$, let $\ok\subset \ol$ be their rings of integers. Let $e$ be
the absolute ramification index of $L$, $s=[\log_p(pe/(p-1))]$.
Let $e_0=[L:K_0]$, where $K_0$ is the maximal unramified
subextension of  $L/K$
 (for ordinary local fields $e_0$ is equal to $e(L/K)$), $l'=s+v_p(e_0)+1$,
 $l=2s+v_p(e_0)+1$.

For an $\ol$-group scheme $H$ we denote the $\ol$-dual of the
module $J/J^2$ by $TH$. Here $J$ is the augmentation ideal of the
coordinate ring of $H$.

\begin{theor}[A]
Let $V$ be an $m$-dimensional Abelian variety over $K$. Suppose
that $V$ has semistable reduction over $L$.

Then $V$ has semistable reduction over $K$ if and only if
  for a certain finite flat  group scheme $H/\ok$ we have
 $TH_\ol\supset (\ol/p^{l}\ol)^m$ (i.e. there exists an embedding);
  besides there exists a monomorphism
 $g:H_K\to\ke[p^{l}]_{V,K}$.

\end{theor}

Note that for $e<p-1$ we have $l=1$; hence Theorem A is a
generalization of Theorem 5.3 of \cite{co} (where only the
potentially good reduction case was considered).

Our technique also allows us to prove the following criteria for
ordinary reduction (see subsection \ref{ordred}) easily.

\begin{theor}[B]
I Let $V$  be an Abelian variety  of dimension $m$ over $K$
 that has good reduction over $L$.
Then the following conditions are equivalent.

 (1) $V$ has good ordinary reduction over $K$

 (2) For a certain  multiplicative type (i.e. dual-\'etale)
group scheme $H/\ok$ we have
 $TH_\ol\approx (\ol/p^{l'}\ol)^m$;
  besides there exists a monomorphism $g:H_K\to \ke[p^{l'}]_{V,K}$.

  (3) For some  $H_K\subset \ke[p^{l'}]_{V,K}$ and $M$ unramified
  over $K$ we have
  $H_M\cong (\mu_{p^{l'},M})^m$. Here $\mu_{p^{l'}}$
is the group scheme of
  $p^{l'}$-th roots of unity.

II Let $V$ be be an Abelian variety  of dimension $m$ over $K$
 that has semistable reduction over $L$.
 Then the following conditions are equivalent.

(1) $V$ has ordinary  reduction over $K$.

(2)  For some multiplicative type group scheme  $H$ over
$\mathfrak{O}_K$
 there exist embeddings of $H_K$ into
$\operatorname{Ker}[p^{l}]_{V,K}$, and of $(\ol/p^l\ol)^m$ into
$TH_\ol$.

(3) For some  $H_K\subset \ke[p^{l}]_{V,K}$ and $M$ unramified
  over $K$ we have
  $H_M\cong (\mu_{p^{l},M})^m$.
\end{theor}

In the first section we remind the definition of the Cartier
modules for formal groups. In this paper we will usually use the
so-called invariant Cartier modules (i.e. $D_F$) defined in
\cite{0} and \cite{01}. We also recall  Oort's definition of the
Cartier module $C(S)$ for a finite local group scheme $S$ (see
\cite{oo}); we remind that the Oort functor is fully faithful. We
remind  the definition of a closed submodule of a Cartier module
and of a separated Cartier module. We also recall connection of
closed submodules with closed subgroups (cf.  \cite{02}). We prove
three technical lemmas on the Cartier modules telated to the
reduction of formal groups and finite group schemes over $\ol$.
The lemmas are used in the proof of the main classification result
(Theorem \ref{mclas}).

Section 2 is dedicated to the main classification theorem. We
describe  the image of the Oort functor; that gives us a full
classification of finite flat local commutative  group schemes
over the rings of integers of complete discrete valuation fields.
The proof that for any $S$ the module $C(S)$ satisfies the
conditions of Theorem  \ref{mclas} is easy. In the proof of the
converse statement we present the given Cartier module as a factor
of $C(F)$ for some formal group $F$; next we prove that a finite
height $F$ can be chosen. The classification immediately implies
that $\ext^1$ for local group schemes coincides with $\ext^1$ for
their Cartier modules.

In section 3 we prove that the tangent space functor can be easily
described in terms of Oort modules. This implies that the minimal
dimension of a finite height formal group $F$ that contains a
given local group scheme $S$ as a closed subgroup is equal to the
minimal number
 of generators for the coordinate ring of $S$.
 The case $m=1$ of this result was essentially proved in \cite{B}
 (though there it was formulated in a  different way).
It was quite important for the program of "taming wild extensions
by means of Hopf  algebras".
 We also prove that a local group scheme of exponent $p^r$ is truncated
  Barsotti-Tate iff $TH\approx (\ol/p^r\ol)^m$.

 In the beginning of  section 4 we remind the main generic fibre result
  and a certain descent result of \cite{02}.
Next  we  prove   certain new descent statements for $p$-divisible
groups. We prove that an Abelian variety has semistable reduction
if and only if its 'formal part' has 'good reduction'. Using these
statements and results on tangent spaces
 we prove a certain  criterion for  good reduction of a potentially good reduction Abelian variety. This criterion gives a positive answer to the question of existence of a finite $p$-adic criterion, which B. Conrad  (see \cite{co}) attributes to N. Katz.
 Two more variants of this criterion were proved in \cite{02}.

In section 5 we prove Theorems A and B.

Theorem A seems to be less convenient than the corresponding
finite $l$-adic criteria (see \cite{sz1}, \cite{sz2}). On the
other hand,  being ordinary is a $p$-adic property. It seems to be
difficult to check it $l$-adicaly; in particular, no $l$-adic
analogue of Theorem B is known.

Note that the usage of tangent spaces makes all formulations (and
certainly the proofs) much different from those in the $l$-adic
case. In particular, we are able to apply the dimension argument
to finite schemes.

One is not able to prove these criteria without using tangent
spaces (even in the case $e<p-1$).

For the convenience of the reader we note that the proofs of
reduction criteria use neither the classification results of
section 2 nor the technical lemmas of subsection \ref{redgreb}.

The author is deeply grateful to prof. Yu. Zarhin for calling his
attention to the case of semistable reduction and to prof. W.
Messing for useful remarks.

\begin{notata}

We keep the notation of the introduction ($K$, $L$, $p$, $e$, $s$,
$e_0$, $l$, $l'$, $\ok$, $\ol$, $V$, $\mu_{p^r}$).

Let $\ovl$ denote the residue field of $\ol$, let $\gm$ denote the
maximal ideal of $\ol$, let $\pi\in \gm$ be some uniformizing
element of $L$.

We introduce $t=v_p(e)+1$.

Let $X=(X_i)=X_1,\dots,X_m$, $x$ be formal variables.

 $F$ will usually denote an $m$-dimensional formal group law over $\ol$,
 $\la_F=(\la_i(X)),\ 1\le i\le m$ denotes the logarithm of $F$;
 $\exp_F$ is the inverse to logarithm (with respect to composition).

In this paper a 'group scheme' will (by default) mean a finite
flat commutative $p$-group scheme (i.e. annihilated by a power of
$p$), $S/\ol$ means a finite flat commutative group scheme over
$\ol$.

 $\car=\car_p(\ol)$ will denote the $p$-Cartier ring over $\ol$
(see below); Cartier module is a module over $\car$.

For a (possibly, non-commutative) ring  $\ga$ we denote by
$M_{m\times n}(\ga)$ the module of $m\times n$-matrices over
$\ga$; $M_m(\ga)=M_{m\times m}(\ga)$.

For finite group schemes $S,T$ we  write $S\nde T$  if $S$ is a
closed subgroup scheme of $T$; we write $S\subset T$ if there is a
group scheme morphism $f:S\to T$ that is injective on the generic
fibre.

We call the number of indecomposable summands of a finitely
generated $\ol$-module $M$ the dimension of $M$.

A formal group $F$ is called a finite height one if $[p]_F$ is an
isogeny, i.e. the scheme-theoretic kernel $\ke[p]_F$ of
$[p]_F:F\to F$ is a finite group scheme. This property is stable
with respect to  base change (see \cite{z}). In particular, $F$ is
finite height over $\ol$ if and only if it is finite height over
some $P\supset \ol$, where the residue field of $P$ is perfect.

For a finite  height formal group $F$ and a group scheme $H$ we
write $H\nde F$ if $H$ is a closed subgroup of $\ke [p^r]_F$ for
some $r>0$ (and hence, a subgroup of the $p$-divisible group of
$F$).

A group scheme or a $p$-divisible group (in particular, a finite
height formal group) is called multiplicative type if its
Cartier-dual is \'etale.

\end{notata}

\section{Cartier modules; some results of previous papers}

In this section we introduce (a modified version of) the classical
Cartier module theory for formal groups; we also recall the result
of Oort and some results of the previous papers. We aso prove a
few technical lemmas.

\subsection{The definition of the Cartier ring and the Cartier modules of formal groups}

For a $\ol$-algebra $Q$ we denote by $\car=\car(Q)$ the ring that
is obtained by factorising $\z\lan\bff, \lan a\ra \ra \lan\lan\bv
\ra \ra$  (here $a\in Q$, we consider non-commutative series over
non-commutative polynomials) modulo the following relations:
\begin{align}
\lan a\ra  \lan b\ra =\lan ab\ra   \text{ for all }a,b\in Q;\
\bff\bv=p\\ \lan a\ra \bv=\bv \lan a^p\ra;\    \lan a^p\ra
\bff=\bff \lan a\ra \text{ for any } a\in Q\\ \lan a\ra +\lan b\ra
=\lan a+b\ra\ +\sum_{n>0} \bv^n \lan r_{p^n}(a,b)\ra   \bff^n,
\label{sum}
\end{align} for all $a,b\in Q$, where $r_{p^n}$ are the polynomials defined in \cite{1}, section 16.2. We need the property $r_{p^n}(0,x)=r_{p^n}(x,0)=0$.

A natural analogue of (\ref{sum}) is valid for $\lan a\ra -\lan
b\ra$.

In \cite{1} $\bv$ was denoted by $V_p$, $\bff$ was denoted by
$F_p$.

Since $\car$ is $\bv$-complete, any finitely generated
$\car$-module is also $\bv$-complete.

If $P$ is a $Q$-algebra we can define the $\car$-module structure
on $P[[\de]]$ in the following way: for $f=\sum_{i\ge 0} c_i\de^i
\in P[[\de]],\ a\in Q$  we define
 $$\bv f=f \de;\ \bff f=\sum_{i>0} pc_i \de^{i-1};\ \lan a\ra f=\sum a^{p^i}c_i\de^i .$$

We also remind that $\car$ is $p$-complete if $Q$ is; if
 $M$ is a $\car(Q)$-module, then $M/\bv M$
has a natural structure of a $Q$-module defined via $$a\cdot
(x\mod\bv M) =\lan a\ra x \mod\bv M $$ for any  $x\in M$.

The main ingredients of the classical Cartier theory were the
modules of curves. For a formal group $F/Q$ and any $i>0$ one can
define the Cartier module structure of $F(x^iQ[[x]])=\ke
(F(Q[[x]])\to F(Q[[x]]/x^i))$. The $p$-typical (abstract) Cartier
module $C(F)$ is a certain direct $\car(Q)$-summand of
$F(xQ[[x]])$. One can describe it explicitly in the case when $Q$
is torsion-free.

For
 $h\in (\mathfrak{A}[[\Delta]])^m$, the coefficients of $h $ are equal to
$h_{i}=\sum_{l\ge 0} c_{il}\Delta^l$, $c_{il}\in\mathfrak{A}$, we
define
\begin{equation}\label{ap}
h(x)=(h_i(x)),\ 1\le i\le m,\text{ where }h_i=\sum_{l>0}c_{il}
x^{p^l}.
\end{equation}

Now let $Q=\ol$. We define $D_F=\{f\in L[[\de]]^m: \exp_F(f(x))\in
\ol[[x]]^m\}$. In particular, for $m=1$ we have $\sum a_i\de^i\in
D_F\iff\exp_F(\sum a_ix^{p^i})\in \ol[[x]]$.

The following statements were proved in  \cite{01} and \cite{0},
also see \cite{1}.

\begin{pr}\label{mcart}
1. $D_F\approx C(F)$ as a $\car$-module.

2. $F\to D_F$ is an embedding of categories of formal groups over
$\ol$ into $\car$-submodules of $L[[\de]]^m$.

Besides, for $\dim F_i=m_i,\ i=1,2$ we have
$$\car(D_{F_1},D_{F_2})=A\in M_{m_2\times m_1}\ol: AD_{F_1}\subset
D_{F_2}.$$

3. If $f:F_1\to F_2$, $f\equiv AX\mod \deg 2$, $A\in M_{m_2\times
m_1}\ol$, then the associated map  $f_*:D_{F_1}\to D_{F_2}$ is the
multiplication by $A$.

4. $D_F\equiv \ol^m\mod\de$.

5. If $P$ is the ring of integers of a complete discrete valuation
field containing $L$ then $D_{F,P}=\car_PD_F\subset PL[[\de]]^m$
and $D_F=D_{F,P}\cap L[[\de]]^m$.

6. The map $D_F\mod\de\to \{f\in xL[[x]]^m:\exp_F(f)\in
\ol[[x]]^m\}/x^2L[[x]]^m$ that is induced by (\ref{ap}) is an
isomorphism.

\end{pr}

The set of all possible $D_F$ can be described by conditions
$D\nde L[[\de]]^m$ (see the definition below), $D\mod\de=\ol^m$
(see  \cite{01}, \cite{02}).

One can easily see (cf. the reasoning in subsection \ref{suflc}
below) that for a large class of $\car$-modules  any $x_i\in M,
1\le i\le m$, such that $\ol (x_i\mod \bv M)=M\mod\bv M$ generate
$M$ over $\car$. Moreover, any $x\in M$ can be presented in the
form $\sum_{1\le i\le m}\sum_{0\le j} \bv^j \lan a_{ij}\ra x_i$.

In particular, this is valid for $M=D_F$. Hence one can choose
$x_i\in D_F,\ x_i\equiv e_i\mod \de,\ 1\le i\le m$
($e_i=(0,\dots,1,\dots,0)$ is the $i$-th basic vector of $\ol^m$).

\subsection{Some results on Cartier modules related to the reduction}\label{redgreb}

The definition of the Cartier ring is functorial. In particular, a
$\car$-module structure on $M$ induces a $\car(\ovl)$-modules
structure on $M/\car\lan \pi\ra  M$.

If $\ovl$ is perfect, the module $C(F)/\car\lan \pi\ra  C(F)$ can
be identified with the usual (covariant) Dieudonne module of the
reduction of $F$.

We denote by $F_\pi$ the formal group law $F_\pi (X,Y)=\pi\ob (\pi
X,\pi Y)$. Its coefficients  obviously belong to $\ol$.

In the proof of the main classification theorem we will need the
following technical results.

\begin{lem}\label{dogreb}
  1. $\car\lan \pi\ra D_F=\pi D_{F_\pi}$.

2. $p^s\pdp\subset \ol[[\de]]^m$.
\end{lem}

\begin{proof}
  1. We can assume that $F$ is $p$-typical (see \cite{1}), i.e. the logarithm $\la=(\la_i)$ of
$F$ satisfies $\la_i=\sum_{1\le j\le m,\ l\ge 0} a_{ijl}
X_j^{p^l}$. Then for $b_i=(b_{ij});\ b_{ij}=\sum a_{ijl}\de^l$ we
have $b_i\in D_F,\ 1\le i\le m$ (see \cite{1}). $b_i\equiv
e_i\mod\de$; hence $b_i$ generate $D_F$ as a $\car$-module.
Similarly $\pi\ob\lan\pi \ra b_i$ generate $D_{F_\pi}$. Since
$\lan\pi \ra b_i \in \lan \pi\ra D_F$, we obtain $\pi
D_{F_\pi}\subset \car\lan \pi\ra D_F$.

On the other hand, $\pi D_{F_\pi}$ is a $\car$-module. We have

$$\lan\pi\ra(\sum_{1\le i\le m}\sum_{0\le j} \bv^j \lan a_{ij}\ra
x_i)= \sum_{1\le i\le m}\sum_{0\le j} \bv^j \lan
a_{ij}\pi^{p^j-1}\ra \lan \pi\ra x_i\in \pi D_{F_\pi}.$$ Hence
$\car\lan \pi\ra D_F\subset \pi D_{F_\pi}$.

2. Easy calculation; see \cite{01}.
\end{proof}

For the same  $x_i$ the following statement is fulfilled.

\begin{lem}\label{greb}
Let $y,z\in D_F$, $z\equiv y\mod \pi D_{F_\pi}$, $y=\sum_{1\le
i\le m}\sum_{0\le j} \bv^j \lan y_{ij}\ra x_i$. Then $z$ can be
presented in the form
\begin{equation}\label{greba}
z=\sum_{1\le i\le m}\sum_{0\le j} \bv^j \lan c_{ij}\ra x_i\text{
where }c_{ij}\equiv y_{ij}\mod\pi.
\end{equation}
\end{lem}
\begin{proof}

We construct inductively (with respect to $l$) a sequence of
$z_l,r_l$ such that $z=z_l+r_l$, whence $z_l$
 satisfy the conditions of (\ref{greba}) with $b_{ijl}$ staying instead of
$c_{ij}$, and $r_l\in \bv^l\pdp$.

We start with $z_0=y$, $r_0=z-y$.

Now suppose that we have $b_l,r_{l,i}$ for some $l\ge 0$.
 We consider the decomposition
$$r_l=\sum_{1\le i\le m}\sum_{0\le j} \bv^j \lan d_{ijl}\pi\ra
x_i.$$ Then according to  (\ref{sum}) we have $$\bv^l (\sum_{1\le
i\le m} (\lan b_{ill}\ra+\lan d_{ill}\pi\ra) x_i)=
\bv^l(\sum_{1\le i\le m}\lan b_{ill}+d_{ill}\pi\ra x_i)+s_l$$ for
certain $s_l\in \bv^{l+1}\pdp$. Hence we may take
$$z_{l+1}=z_l+\sum_{1\le i\le m}(\lan b_{ill}+d_{ill}\pi\ra-\lan
b_{ill}\ra) x_i.$$ Since $D_F$ is $\bv$-complete, passing to the
limit yields the assertion.
\end{proof}

We denote by $W$ the twisted power series ring
$\ovl[[\bv']]=\sum_{i\ge 0}\bv^ix_i$; the multiplication is
defined by the rule $x\bv =\bv x^p,\ x\in \ovl$.

Let $\mw$ be a module over $\car'=\car(\ovl)$ (hence also a
$\car$-module). Since the characteristic of $\ovl$ is $p$, we have
$\bv\bff=p=\bff\bv$ in $\car'$ (see \cite{1}).

Suppose that any $x\in\mw$ can be presented in the form
$\sum_{1\le i\le m}\sum_{0\le j} \bv^j \lan y_{ij}\ra x_i$ for
$y_{ij}\in\ovl$. For any $z=(z_i)\in W^m$, where $z_i=\sum \bv^j
z_{ij}$,
 we consider $$(z)(x)=\sum_i\sum_{j\ge 0} \bv^j\lan z_{ij}\ra x_i.$$
Note that $(z_1+z_2)(x)\neq (z_1)(x)+(z_2)(x)$.

We introduce a topology on $W$ whose basic neighbourhoods of $0$
are $V^iW$. Obviously, if $s_i\to 0$ in $W^m$ then $\sum (s_i)(x)$
converges in $\mw$.

\begin{lem}\label{grebanaya}
Let $X$ be a $W$-submodule of $W^m$; we consider the set
$Y=(X)(x)=\{(z)(x);\ z\in X\}$. Suppose that $\bff x_i\in Y$ for
$1\le i\le m$.

Then the following statements are fulfilled.

1. For any $y,z\in X$  there exist $u,u'\in  X$ such that
$(u)(x)=(y)(x)+(z)(x)$, $(u')(x)=(y)(x)-(z)(x)$.

Besides, if $z\in\bv^r W^m$ for $r>0$ then we can choose
$u,u'\equiv y\mod \bv^r W^m$.

2. $Y$ is a $\car'$-submodule of $\mw$.

3. For any $y\in X$, $z\in W^m\setminus X$ there exists a $u\in
W^m\setminus Y$ such that $(u)(x)=(z)(x)-(y)(x)$.

\end{lem}
\begin{proof}
1. The proof is similar to the  proof of the previous Lemma.

First we construct $u$.

We denote $(y)(x)+(z)(x)$ by $a$.

We construct  inductively (with respect to $l$)  a sequence of
$b_l,r_{li}\in X,\ l,i\ge 0$, such that $\lim_i r_{l,i}=0$ for any
$i\ge 0$, $(b_l)(x)+\sum_i (r_{l,i})(x)=a$, and $r_{l,i}\in
\bv^lX$.

We start with  $b_0=y$, $r_{0,0}=z$, $r_{0,i}=0$ for $i>0$.

Now suppose that we have $b_l,r_{l,i}$ for some $l\ge 0$.

For $a\in \olv,\ c\in W^m$ we have $\lan a\ra(c)(x)=(ac)(x)$.
Hence \begin{equation}\label{eqy} \lan a\ra Y\subset Y\text{ for
any } a\in\ovl\end{equation}

Now arguing similarly to the  proof of Lemma \ref{greb}, we obtain
that we can take $b_{l+1}=b_l+\sum_{i\ge 0} r_{l,i}$. Indeed,
according to (\ref{sum}) $$(b_{l+1})(x)-(b_{l})(x)-\sum_{i\ge 0}
(r_{l,i})(x)$$ can be expanded as a sum of   $(w_{l+1,i})(x)$ for
some $w_{l+1,i}\in \bv^{j_{li}}\lan\ovl\ra \bff W^m\subset
\bv^{j_{li}}X$. Here each $j_{li}$ is greater than $l$, $\lim_i
j_{li}=\infty$. Hence we can take $r_{l+1,i}=w_{l+1,i}$.

Lastly we take $u=\lim_l b_l$.

Obviously, for $z\in\bv^r W^m$ this construction gives us a $u$
that is  congruent to $y$ modulo $\bv^r W^m$.

The construction  of $u'$ is obtained from those for $u$ by
replacing some $+$ signs by $-$.

2. We have $\bv Y\subset Y$. Part 1 states that $Y\pm Y\subset Y$.
We also have (\ref{eqy}). It remains to check that $\bff Y\subset
Y$.

Let $d=\sum \bv^ld_l\in W^m;\ d_l\in \ovl^m$. Then we have $$\bff
(d)(x)=\sum_{l\ge 0} \bv^{l} (\bff d_l)(x).$$ Besides, $\bff
(d_l)(x)\in Y$. According to the assertion of part 1, there exists
a sequence $u_i\in X,\ i\ge 0$ such that $(u_i)(x)=\sum_{0\le j\le
i} \bv^{j} (\bff d_j)(x)$, and $u_{i+1}\equiv u_i\mod \bv^{i+1}$.
Hence there exists a limit  $u=\lim u_i$; we have $(u)(x)=\bff
((d)(x))$.

Thus $\bff Y\subset \bff\mw\subset Y$.

3. The proof is very similar to those of of part 1. We construct
 inductively (with respect to $l$)
$b_l,r_{li}\in X,\ l,i\ge 0$, such that $(b_l)(x)+\sum_i
(r_{l,i})(x)=(z)(x)-(y)(x)$; $\lim_i r_{l,i}=0$, and $r_{l,i}\in
\bv^lX$.

In order to pass from $b_l\notin X$ to $b_{l+1}\notin X$ we use
the fact that $(W^m\setminus X)\pm X=W^m\setminus X$.
\end{proof}

\subsection{Oort modules of group schemes}

Let $S$ be a local group scheme over $\ol$; let $0\to S\to F\to
G\to 0$ be its resolution by means of finite height formal groups.
We define $D(F)=\cok(C(F)\to C(G))$.

In the paper \cite{oo} the following result was proved.

\begin{pr}
 $S\to C(S)$ is a well defined functor on the category of local (finite flat commutative) group schemes over $\ol$; it defines an embedding of this category into the category of $\car$-modules.
\end{pr}

We call $C(S)$ the Oort module of $S$.

\subsection{Closed submodules and closed subschemes}\label{clo}

The following definition played a crucial role in \cite{02}.

\begin{defi}
1. For $\car$-modules $M\subset N$ we write $M\nde N$, if for any
$x\in M, \bv x\in N$ we have $x\in  M$. We call $M$ a {\bf closed}
submodule of  $N$.

2. $\car$-module $N$ is called {\bf separated} if $\ns\nde N$,
i.e. $N$ has no $\bv$-torsion.
\end{defi}

For Cartier modules (i.e. modules over $\car$) one can easily
prove the following facts (cf. \cite{02}).

\begin{pr}
\label{mclo}

I1. If $M_i\nde M$ for $i\in I$, then $\bigcap_{i\in I} M_i\nde
M$.

2. If $M_2\nde M_1$, $M_1\nde M$, then $M_2\nde M$.

3. If $M_1\nde M$,
 $M_2\subset M_1$, then $ M_1/M_2\nde M/M_2$.

II If $f:N\to O$ is a $\car$-morphism, $M\nde O$, then
$f\ob(M)\nde N$.

\end{pr}

For a Cartier module $N$,  $O\subset N$, we denote
 by $\cl_N(O)$ the smallest closed Cartier-submodule of $N$ that contains $O$.
 In particular, $\cl_N(O)\supset \car O$.

 Obviously, if $M\subset N$ are $\car$-modules, then $\cl_N(M)=M\iff \cl_{N/M}\ns=N/M$.

Now we remind  the connection of closed modules and closed
subgroup schemes.

\begin{theo}\label{clsc}
I For a local group scheme $S$ the module $C(S)$ is separated.

II1. Closed submodules of $C(S)$ are in one-to-one correspondence
with closed subgroup schemes of $S$.

2. If $M=C(S)$, $N\nde M=C(H)$, where $H\nde S$, then $M/N\approx
C(S/H)$.

3. Conversely, exact sequence (as fppf-sheaves, i.e. the inclusion
is a closed embedding) of local schemes induce exact sequences of
Oort modules.

III If $f:S\to T$ is a local scheme morphism, then $\ke f_*=C(\ke
f)$, where $f_*$ is the induced Oort module homomorphism; we
consider the kernel in the category of flat group schemes.

IV Let $M$ be a closed finite index submodule of $C(F)$ where $F$
is a finite height formal group. Then $M$ is isomorphic to $C(G)$,
where $G$ is a formal group that is isogenous to $F$.

V If $M\nde L[[\de]]^m$, $D_F\subset M$, and $M\mod\de=\ol^m$,
then $M=D_F$.
\end{theo}

Parts I-III were proved in \cite{02}, IV and V were proved in
\cite{01} (in  slightly different formulations).

\section{The main classification theorem; extensions of group schemes}

This section is dedicated to classification of group schemes in
terms of their Oort modules. We describe the image of the
Cartier-Oort functor, thus giving a complete classification.

\subsection{Formulation}

\begin{theo}\label{mclas}
I A $\car$-module $M$ is isomorphic to $C(S)$ for $S$ being a
finite flat commutative local group scheme over $\ol$ if and only
if $M$ satisfies the following conditions.

1. $M/\bv M$ is a finite length $\ol$-module.

2. $M$ is separated.

3. $\cap_{i\ge 0} \bv^i M=\{0\}$.

4. $M=\cl_M(\lan\pi\ra M)$.

II The minimal dimension of a finite height formal group $F$ such
that $S$ can be embedded into $F$ is equal to $\dim_\ol(C(S)/\bv
C(S))$.

III $M=C(\ke [p^r]_F)$ for an $m$-dimensional formal group $F$ if
and only if in addition to the conditions of I, $p^r M=0$ and
$M/\bv M\approx (\ol/p^r\ol)^m$.
\end{theo}

We will denote $\dim_\ol(M/\bv M)$ by $\dim M$.

\subsection{The proof of necessity in Theorem \ref{mclas}}\label{nec}

Let $M=C(S)$.

Suppose that $S=\ke f$, where $f:F\to G$ is an isogeny of finite
height formal group laws of dimension $m$. Then for $M=\cok f_*$
we have $M/\bv M=D_G/(\de D_G+AD_F)$ where $f\equiv AX\mod \deg
2$. Since $D_F\nde L[[\de]]^m$, we have
\begin{equation}\label{as}\begin{aligned}
M/\bv M=D_G/(\de D_G+AD_F)\\ =(D_G\mod\de)/(AD_F\mod\de)\approx
\ol^m/A\ol^m.
\end{aligned}
\end{equation}
 We obtain that $M$ satisfies I1. We also see that the minimal dimension of a finite
 height formal group $F$ such that $S$ can be embedded into $F$ is not
 less than $\dim(C(S)/\bv C(S))$. Besides, if $S=\ke [p^r]_F$,
then $C(S)/\bv C(S)\approx (\ol/p^r\ol)^m$.

I2 for $M$ is part I of Theorem \ref{clsc}.

Since $f_*D_G\nde D_F$ we obtain $$C(S)/\bv^r C(S)\cong D_G/(\bv^r
D_G+ f_*D_F)\cong (D_G/\bv^r D_G)/ (f_*D_F/\bv^r f_*D_F).$$ Hence
$$\varprojlim M/\bv^rM\cong \prli (D_G/\bv^r D_G)/ (f_*D_F/\bv^r
f_*D_F)=D_G/f_*D_F\cong M.$$ We immediately obtain I3 for $M$.

Now we prove I4. Suppose that $N=\cl_M(\lan\pi\ra M)\neq M$. Then
$N$  corresponds to some $T\nde S$. We obtain that the module
$U=\{x\in D_G:\ x\mod f_* D_F\in N\}\supset \pi D_{G_\pi}$. Since
$p^rS=0$ for some $r>0$, we have $p^rM=0$; hence $U$ is a finite
index submodule of $D_G$. Therefore $U$ corresponds to a formal
group that is isogenous to $G$ (see Theorem \ref{clsc}, part IV).
Since these conditions are formulated in terms of  formal groups
and their Cartier modules, they are not affected by complete
discrete valuation  fields extensions (see Proposition
\ref{mcart}, part 5). Hence we may assume that $\ovl$ is perfect.
We obtain that $M/\lan\pi\ra M$ is the (covariant) Dieudonne
module of the reduction of $S$. Since $S$ is local, its reduction
is local also. It is well known that in this case $M/\lan\pi\ra M$
is $\bv$-torsion, which contradicts the existence of $T$.

\subsection{The proof of sufficiency in Theorem \ref{mclas}:
 construction of a certain formal group}\label{suflc}

Suppose that $M$ satisfies conditions  of Theorem \ref{mclas} and
$\dim M=m$.
 We first prove that there exists a formal group $F$, $\dim F=m$,
 such that $M\approx D_F/C$, where $C\nde D_F$.
 Next we verify that one can choose $F$ being a finite height formal group.

We choose representatives $l_1,\dots,l_m\in M$ for (some)
$\ol$-generators of $M/\bv M$. We check that any $x\in M$ can be
expressed as a sum
\begin{equation}\label{csum}
\sum_{1\le i\le l}\sum_{0\le j} \bv^j \lan a_{ij}\ra l_i.
\end{equation}

Indeed, for any $x\in M$ there exist  $a_i\in \ol$ such that $x=
\sum_{0\le i}  \lan a_{i}\ra l_i+Vx_1$ for some $x_1\in M$.
Repeating this procedure $r$ times we obtain that $$x=\sum_{1\le
i\le l}\sum_{0\le j<r} \bv^j \lan a_{ij}\ra l_i+V^rx_r$$ for some
$x_r\in M$. Since $\sum \bv^j \lan a_{ij}\ra$ converges in $\car$
for any $a_{ij}\in \ol$, we obtain that there exists an $y\in M$
that can be presented in the form (\ref{csum}) and $x-y\in V^r M$
for any $r>0$. Now condition I3 for $M$ implies that $x=y$.

Next we present $\bff l_i$ in the form (\ref{csum}) for all $1\le
i\le m$. We obtain
\begin{equation}\label{presfr}
\bff l_i=\sum_{1\le k\le l}\sum_{0\le j} \bv^j \lan a_{ijk}\ra
l_k.
\end{equation}
According to section 27.7 of \cite{1}, there exists an
$m$-dimensional formal group law $F$ such that $$C(F)\approx
(\bigoplus \car e_i)/\sum \car (\bff e_i-\sum_{1\le k\le
l}\sum_{0\le j} \bv^j \lan a_{ijk}\ra e_k)$$ (here $e_i$ denote
certain $\car$-generating elements). We obtain that there exists a
$\car$-homomorphism $h$ of $C(F)$ to $M$ that maps $e_i$ to $l_i$.
Since $l_i$ generate $M$ over $\car$,  $h$ is onto.

Since $M$ is separated, $\ke h\nde C(F)$.

\subsection{Choosing a finite height formal group}

Now we prove that we can choose $F$ being a finite height group.
Since the height is determined by the reduction of $F$, it is
sufficient to consider the residues of $a_{ijl}$, which we denote
by $\overline{a_{ijk}}$.
 Again we denote
$M/\car \lan \pi\ra M$ by $\mw$. According to Lemma \ref{greb},
the calculation of $\overline{a_{ijk}}$ can be done in $\mw$. We
denote the images of $l_i$ in $\mw$ by $\olv_i$.

We embed $\ovl$ into a perfect field $U$. It is well known that
$F$ has finite height if and only if the matrix $B=b_{ik}$, where
$b_{ki}=\sum_{j\ge 0}\overline{a_{ijk}}\bv^{j}$ is non-degenerate
over the skew-field $W'=U((V))'=\sum_{i\gg -\infty}\bv^ic_i$. Here
$W'\supset W$ is the skew-field of twisted Laurent series over
$U$, the multiplication is defined by the rule $u^p\bv =\bv u,\
u\in U$.

As in subsection \ref{redgreb}, for any $z=(z_i)\in W^m$ we define
$(z)(\olv)=\sum z_i\olv_i\in \mw$.

 Let $v\le m$ be the maximal possible rank of $B$ for all possible choices of the expansion (\ref{presfr}).

Suppose that $v<m$. We fix the numbers $i_1,\dots,i_v$ and the
corresponding columns $b_{i_l}$ of $B$ such that $b_{i_l}$ are
$W'$-independent. We obtain that all other columns belong to the
$W'$-span of $b_{i_l}$ for any possible choice of their
coefficients. We denote $(\sum W'b_{i_l})\cap W^m$ by $X$; denote
$(X)(\olv_i)\subset \mw$ by $Y$. Since $W'(W^m)\neq \sum
W'b_{i_l}$, we have $X\neq W^m$. We have $\bv X\subset X;\ \bv
(W^m\setminus X)\subset W^m\setminus X$. Since (for some choice of
$F$) all columns of $B$ belong to $X$, we have $\bff \olv_i\in Y$
for $1\le i\le m$. Hence we can apply Lemma \ref{grebanaya} with
$\olv_i$ staying for $x_i$.

We choose a column $b=b_h$ for $h\notin \{i_l\}$. We have $b\in
X$. We prove that we can modify $F$ so that the new $b$ will not
belong to $X$.

We choose certain $z\in W^m\setminus X$; the choice depends on
whether $Y=\mw$.

In the case $Y=\mw$ we choose any $z\in W^m\setminus X$.

Suppose that $Y\neq\mw$. According to part 3 of Lemma
\ref{grebanaya},  $Y$ is a $\car'$-module. Since the closure of
$\ns$ in $\mw$ is equal to $\mw$, we obtain that $Y\not\nde \mw$.
Hence there exists $y\in \mw\setminus Y$ such that $\bv y\in \mw$.
We expand $y=(z)(\olv_i)$, $z=(z_i),\ z_i\in W$. We have $z\notin
X$.

In both cases we obtain that there exists a $c\in X$ such that
$(c)(\olv)=\bv (z)(\olv)$. According to part 3 of Lemma
\ref{grebanaya} there exists a $g\in W^m\setminus X$ such that
$(g)(\olv)=(c)(\olv)-\bv (z)(\olv)=0$. According to part 1 of
Lemma \ref{grebanaya}, there exists a $d\in X$ such that
$(d)(\olv)=-(b)(\olv)$. Applying part 3 of Lemma \ref{grebanaya}
for $z=g$, $y=d$, we obtain that $b$ can be replaced by $b'\notin
X$.

Hence the case $v<m$ is impossible. Therefore $v=m$. Thus we can
choose $F$ so that $B$ is invertible over $W'$, hence $F$ is a
finite height formal group.

\begin{rema}
1. In the perfect residue field case the reasoning could be
simplified by using the fact that $\bv^u M\subset \car \lan \pi
\ra M$ for some $u>0$. This is no longer true in the imperfect
residue field case.

2. Another way of proof is possible. Suppose that for a map
$f:F\to G$ of (not necessarily finite height) formal groups we
have $C(G)/f_*C(F)\approx M$ for $M$ satisfying the conditions I1
-- I4. Then one can  note that $M$ can be obtained from $\car
\lan\pi\ra M$ by means of successive adjoining of roots of
equations of the type $\bv y=x$. Hence $C(\overline{G})$ can be
obtained from the  Cartier module for the reduction of $f_*C(F)$
by successive adjoining the roots of the equation of the type
$\bv(x_i)=(y_i)$, $1\le i\le m$. Next the Oort type reasoning
could be applied to prove that $\ke f$ is a finite flat group
scheme and $M=C(\ke f)$.
\end{rema}

\subsection{The end of the proof of parts I and II}

We take any  $m$-dimensional finite height $F$ such that $M$ can
be presented as a factor of $D_F$ (as a $\car$-module).

We prove that the exponent of $M$ is finite, i.e. for some $u>0$
we have $p^uM=0$. This is equivalent to $p^uD_F\subset \ke h$,
where $h:D_F\to M$  is the map constructed in subsection
\ref{suflc}.

First we check that for some $r>0$ we have $p^r\mw=0$. We recall
that $$\mw=\sum_{1\le i\le m}\car' \olv_i=\cl_{\mw}\ns.$$ Hence
there exists a sequence of $\car'$-modules $\ns=M_0\subset
M_1\subset M_2\subset\dots\subset M_r=\mw$ such that for any $0\le
i< r$ we have $M_{i+1}=M_i+\car x_i$ for some $x_i\in M_{i+1}$
satisfying $\bv x_i\in M_i$. We have $px_i=\bff\bv x_i\in M_i$.
Therefore $pM_{i+1}\subset M_i$. Thus $p^r\mw=0$.

Next we verify that $p^g\pdp \subset\ke h$ for some $g>0$. Since
$M/\bv M$ is a finite length $\ol$-module, we obtain that $\ke
h\mod\de$ is a finite index submodule of $D_F\mod\de=\ol^m$. If
$x=\sum_{i\ge 0} x_i\de^i,\ x_i\in L^m$, then $\bff\bv x-\bv\bff
x=p x_0$. Hence for some $d>0$ we have $p^d\ol^m\subset \ke h$;
therefore $p^d\ol[[\de]]^m\subset \ke h$. Using part 2 of Lemma
\ref{dogreb} we obtain $p^{d+s}\pdp \subset\ke h$.

We have an exact sequence $0\to \pdp/(\pdp \cap\ke h)\to
M\to\mw\to 0$. Thus we can take $u=r+g$.

Since $F$ is finite height,  the module $D_F/p^u D_F$ corresponds
to a finite flat local group scheme $\ke [p^u]_F$.

Since $\ke h\nde D_F$, we obtain $\ke h/p^u D_F\nde D_F/p^u D_F$.
Hence, according to part II2 of Theorem \ref{clsc}, $M$
corresponds to some group scheme $S$ which is a factor of $\ke
[p^u]_F$ by a closed subscheme.

\subsection{The proof of sufficiency in Theorem \ref{mclas}, part III}

Again we choose any  $m$-dimensional finite height $F$ such that
$M$ can be presented as a factor of $D_F$.

It is sufficient to prove that for $F,h$  constructed above we
have $p^r D_F= \ke h$. Since $p^r M=0$, we have $p^r D_F\subset
\ke h$.

Since $M/\bv M\approx (\ol/p^r\ol)^m$, using (\ref{as}) we obtain
$\ke h\mod\de =(p^r\ol)^m$. Hence we have $p^r D_F\subset \ke h$,
$p^r D_F\mod\de=\ke h\mod\de$. According to part V of Theorem
\ref{clsc}, we have $p^{-r}\ke h= D_F$. Hence $p^r D_F= \ke h$.

\subsection{Extensions of group schemes}

Theorem \ref{mclas} easily implies the following result for
$\ext^1$ in the category of group schemes.

\begin{pr} If $S,T$ are local, then $\ext^1(S,T)=\ext^1_{\car}(C(S), C(T))$.
Here we consider extensions in the category of finite flat group
schemes, whence the definition of an exact sequence is the same as
always in this paper.
\end{pr}
\begin{proof}
Since $S,T$ are local, any extension of $T$ by $S$ is also local.

We have to prove that conditions of part I of Theorem \ref{mclas}
are preserved by extensions of $\car$-modules.

 Suppose that we have an exact sequence $0\to C(T)\to M\to C(S)\to 0$ of Cartier-modules.
We check that $M$ satisfies conditions I1 -- I4.

 The sequence
$C(T)/\bv C(T)\to M/\bv M\to C(S)/\bv C(S)\to 0$ is right-exact
and we obtain I1.

Since $C(S)$ is separated, we obtain $C(T)\nde M$. Hence $\ns\nde
C(T)\nde M$ and we obtain I2.

Since $\cap_{i\ge 0} \bv^i C(S)=\{0\}$,   for any $x\in \cap_{i\ge
0} \bv^i M$ we have $x\in C(T)$. Since $C(T)\nde M$, we also
obtain $x\in \cap_{i\ge 0} \bv^i C(T)$. Now condition I3 for
$C(T)$ implies $x=0$. We obtain I3 for $M$.

Since $C(T)\subset M$, we obtain $\lan\pi\ra C(T)\subset
\lan\pi\ra M$. Hence $\cl_M(\lan\pi\ra M)\supset C(T)\nde M$.
Hence $$\cl_M(\lan\pi\ra M)=\{x\in M:\ x\mod C(T)\in \cl_C(S)(\car
\lan\pi\ra C(S))\}=M.$$ We obtain I4 for $M$.

\end{proof}

\section{The tangent space of a finite group scheme}

\subsection{Definition;  expression in terms of the Oort module}
We introduce a natural definition of the tangent space $TS$ for a
finite group scheme $S$.

\begin{defi}\label{tgs}
For a finite flat group scheme $S$ we denote by $TS$ the
$\ol$-dual of $J/J^2$ (i.e. $\homm_\ol (J/J^2,L/\ol)$), where $J$
is the augmentation ideal of the coordinate ring of $S$.
\end{defi}

\begin{rema}\label{rts}
1. Obviously, the tangent space is an additive functor on the
category of finite group schemes.

2. The definition implies that the (suitably defined) tangent
space of $S_P=S\times_{\spe \ol} \spe P$ is equal to $TS\otimes_
\ol P$. Here $P$ is any (unitial commutative) $\ol$-algebra.

3. It is well known that the tangent space of a group scheme is
equal (i.e. naturally isomorphic) to the tangent space of its
local part.
\end{rema}

\begin{pr}\label{dts}
1. $TS$ is naturally isomorphic to $C(S_0)/\bv C(S_0)$, where
$S_0$ is the local part of $S$.

2. A local group scheme homomorphism $f:S\to T$ is a closed
embedding  if and only if the induced map of the tangent spaces is
an embedding.

3. If $0\to H\to S\to T\to 0$ is an exact sequence of local group
schemes (in the category of fppf-sheaves, i.e. $H\nde S$) then the
corresponding sequence of tangent spaces is also exact.

4. If $H\subset S$ then the $\ol$-length of $TH$ is not greater
than the length of $TS$; they are equal iff $H=S$.
\end{pr}

\begin{proof} 1. Since the local part is functorial, it is sufficient
to prove the statement for $S_0=S$. We denote the isomorphism we
want to construct by $i_S$.

For a formal group law (i.e. for a formal Lie group) $F$ let $J_F$
denote the augmentation ideal of the coordinate ring of $F$. It is
well known that $J_F/J_F^2$ is the cotangent space of $F$ (at $0$)
i.e. we have an isomorphism
 $$TF= \homm_\ol (J_F/J_F^2,\ol)\cong TF'=\ke (F(x\ol[[x]])\to F(x\ol[[x]]/x^2))$$ that is canonic and functorial.
Via $\log_F$ the module $TF'$ is canonically isomorphic to $\{f\in
xL[[x]]^m:\exp_F(f)\in \ol[[x]]^m\}/x^2L[[x]]^m$. Now applying
part 6 of Proposition \ref{mcart} we obtain $T_F\cong D_F\mod\de$.
 Hence we have a functorial isomorphism $i_F: D_F/\bv D_F\to TF$.

Let $0\to S\to F\stackrel{h}{\to} G\to 0$ be a resolution of $S$
via finite height formal groups of dimension $r>0$. Since $h_*
D_F\nde L[[\de]]^r$, we have $h_* D_F\nde D_G$. According to
(\ref{as}) we have $$C(S)/\bv C(S)=D_G/(h_* D_F+\bv D_G)=(D_G/\bv
D_G)/h_* (D_F/\bv D_F).$$ Besides $h_*$ is injective on $(D_F/\bv
D_F)$. The standard schematic construction of the kernel of a
formal group scheme homomorphism gives us $J\cong J_F/h^*J_G$.
Moreover,  $h^* J_G^2=h^* J_G\cap J_F$. Therefore we can define
$i_S$  by means of $i_G$.

We note that the construction of $i_S$ is functorial with respect
to  morphisms of isogenies of formal groups (i.e. with respect to
commutative squares of morphisms of formal groups).

It remains to check that $i_S$ does not depend on the choice of
resolution for $S$ and is functorial with respect to group scheme
morphisms.

Let $0\to S\to F_1\stackrel{h_1}{\to} G_1\to 0$ and $0\to S\to
F_2\stackrel{h_2}{\to} G_2\to 0$ be two different resolutions of
$S$ by means of formal groups.

It was proved in \cite{oo} that there exists a resolution $0\to
S\to F\stackrel{h}{\to} G\to 0$ such that for $i=1,2$ the
following commutative diagram  can be constructed:
\begin{equation}\label{dia}
\begin{CD}
S@>{}>>F_i@>h_i>> G_i\\ @VV{\id_S}V@VV{}V @VV{}V\\ S@>{}>>F@ >h>>
G
\end{CD}
\end{equation}
the rows are exact.

Hence $i_S$ for first two resolutions are the same as the one for
the third resolution.

Now let $f:S\to T$ be a morphism of (finite flat local
commutative) group schemes. Then according to \cite{oo}, we can
construct a commutative diagram

\begin{equation}\label{bdia}
\begin{CD}
S@>{}>> F_1@>{}>>G_1\\ @VV{f}V @VV{}V @VV{}V\\ T@>{}>>
F_2@>{}>>G_2
\end{CD}
\end{equation}
where the rows are resolutions of $S$ and $T$. Since $i_S$ is
functorial with respect to  morphisms of isogenies of formal
groups, we obtain that $i_S$ is natural.

2. If $H\nde S$ then $C(H)\nde C(S)$, hence we obtain that $
C(H)\cap \bv C(S)=\bv C(H)$.
 Therefore the kernel of the map $C(H)/\bv C(H)\to C(S)/\bv C(S)$ is zero.

We prove the inverse implication.

We may decompose  $f$ as $g\circ i\circ h$, where $g$ is a closed
embedding that corresponds to the injection $\imm f_L\to S_L$, $i$
is bijective on the generic fibre, and $h$ is epimorphic  (see
Proposition \ref{ray} below).

Since the induced tangent space map $g_*$ is injective,  we may
assume that $f_L$ is surjective (i.e. $g=\id_S$).

If the composite $i_*\circ h_*$ of the maps of tangent spaces is
injective then
 $h_*$ is injective.

If we have an exact sequence $0\to B\to H\to S\to 0$, then
$\ns\neq TB\subset TH$; besides $TB$ maps to $0$ via the map
$TH\to TS$. Hence $h$ is an isomorphism.

It remains to consider the case $H\subset S$, $H_L=S_L$. The map
$TH\to TS$ is injective iff $ C(H)\cap \bv C(S)=\bv C(H)$. Then
for $x\in C(H),\ \bv x\in C(S)$ we obtain $\bv x= \bv y,\  y\in
C(H)$. Since $C(S)$ is separated, we obtain $x=y\in C(H)$. Hence
$C(H)\nde C(S)$. Applying Theorem \ref{clsc} we obtain $H\nde S$.

3. From Theorem \ref{clsc} we immediately obtain that the sequence
$TH\to TS\to TT\to 0$ is right-exact. It remains to apply part 2.

4. We can assume that $H$ and $S$ are local.

Suppose that we have exact sequences $0\to H_0\to H\to H_1\to 0$
and $0\to S_0\to S\to S_1\to 0$, where $H_i\subset S_i$ for
$i=0,1$. Then, according to the assertion of part 3, we obtain
that it is sufficient to check the assumption for the pairs
$H_i\subset S_i$. Applying the results of Raynaud (see Proposition
\ref{ray} below) we reduce the statement to the case $H_L=S_L$,
$H_L$ is $L$-irreducible. Since the assumption is not affected by
finite extensions of $L$ (see Remark \ref{rts}), we may also
assume that the order of $H$ is $p$. Then the coordinate ring of
can $H$ be presented as $\ol[x]/x^p-ax$ for some $a\in\ol$  (see
the classification in \cite{ot}), the augmentation ideal is equal
to $(x)$. Then $TH\approx \ol/a\ol$. It remains to note that $\spe
(\ol[x]/x^p-ax)\subset \spe\ol[x]/x^p-bx $ implies $a\mid b$;  if
$b\sim a$ then the inclusion is an isomorphism.
\end{proof}

\begin{rema}
One can also prove the equality $d=l{\#S}$ for $l$ being the
length of $TS$, $d$ being the valuation of the discriminant of the
coordinate ring of $S$. This generalizes Tate's formula for
$d(\ke[p^r]_F)$ (see \cite{ta}).
\end{rema}

\subsection{The dimension of a group scheme}

Now we rewrite parts II and III of Theorem \ref{mclas} in terms of
tangent spaces.

\begin{theo}\label{dimgs}
I For a local group scheme $S$ the following numbers are equal.

1. The $\ol$-dimension of $J/J^2$.

2. The $\ol$-dimension of $C(S)/\bv C(S)$.

3. The minimal dimension of a finite height formal group $F$ such
that $S\nde F$.

II $S$ is equal to $\ke [p^r]_F$ for some $m$-dimensional finite
height formal group $F/\ol$ if and only if $p^r S=0$ and
$TS\approx (\ol/p^r\ol)^m$.
\end{theo}
\begin{proof} Immediate by combining Theorem \ref{mclas} and Proposition \ref{tgs}.
\end{proof}

It seems natural to call the number $\dim_{\ol} (J/J^2)$ the
dimension of $S$.

\begin{rema}
1. We obtain that the minimal number of generators for the
coordinate ring of $S$ is equal to the minimal possible dimension
of $F$, $S\nde F$. Obviously, the dimension of $F$ cannot be
smaller than the minimal number of generators; yet the inverse
inequality is much more difficult and was not known previously.

The case $m=1$ of this result was proved in \cite{B} (though there
it was formulated in a different way).

2. It seems very probable that the equality I1=I3 can be
generalized to not necessarily local $p$-group schemes (and
$p$-divisible groups). In order to do this one should modify the
Cartier module theory so that it would describe not necessarily
local group schemes. This would probably be done in one of the
following papers.
\end{rema}

\section{Descent; a finite  good reduction criterion}

In the paper \cite{02} it was proved that one can check whether an
Abelian variety of potentially good reduction has good reduction
over $K$ knowing $\ke [p^{l'}]_{V,K}$. In this section we modify
that criterion using the results of the previous section. We also
prove that one can reduce the study of semistable reduction for an
Abelian variety to the study of its 'formal part' (see subsection
\ref{pdual} below).

\subsection{The generic fibre results: reminder}

We often use the following fact (see \cite{re}).

\begin{pr}\label{ray}
 Closed subgroup schemes of $S/\ol$ are in one-to-one correspondence with closed subgroup schemes of $S_L$.
\end{pr}

In particular, for any set of  $S_i\nde S$ there exists a
smallest (in the sense of  $\subset$) group scheme $H$ such that
$S_i\nde H$. $H$ is also a closed subgroup scheme of $S$.

Now we remind the main result of \cite{02}. It is a finite
analogue of fullness of the generic fibre functor for
$p$-divisible groups (proved by Tate). Besides, it implies Tate's
result (see \cite{ta}) immediately.

\begin{theo}\label{maingf}
If $S,T$ are $\ol$-group schemes, and $g:S_L\to T_L$ is a morphism
of $L$-group schemes, then there exists an $h:S\to T$ over $\ol$
such that $h_L=p^s g$.

\end{theo}

Certainly, $h$ is unique. One easily checks that the result is
sharp, i.e. the value of $s$ its the smallest one possible.

For $e=e_L<p-1$ Theorem \ref{maingf} implies that the generic
fibre functor on the category of (finite flat commutative) group
schemes is full. This is the central result of \cite{re};
Raynaud's methods cannot be used in the case $e\ge p-1$.

\subsection{Descent lemmas: reminder}

In \cite{02} the following result  was proved by means of  an
explicit descent reasoning for $D_F$.

\begin{pr}\label{oldde}
 Let $F$ be a finite height formal group over $\ol$. Suppose that
its generic fibre
 (as a  $p$-divisible group) $F_L=F\times_{\spe \ol} \spe L$ is defined over $K$,
i.e. there exists a $p$-divisible group $Z_K$ over $K$ such that
\begin{equation}\label{cong}
Z_K\times_{\spe K} \spe L \cong F_L.
\end{equation}
Suppose that  $\ke [p^t]_Z\approx T\times_{\spe \ok} \spe K$
  for some group scheme $T/\ok$;
this isomorphism combined with the isomorphism  (\ref{cong}) is
the generic fibre of a certain isomorphism
 $T \times_{\spe \ok} \spe \ol \cong \ke [p^t]_{F}$.
Then $Z_K\approx F'_K$ for some formal group $F'/\ok$.

\end{pr}

\subsection{Descent for $p$-divisible groups in terms of tangent spaces}

We remind that (by definition) the dimension of a $p$-divisible
group   over $\ol$ is the dimension of its local part (as a formal
group law).

\begin{pr}\label{degr}
Let $V$ be a $p$-divisible group  over $K$, let $Y$ be a
$p$-divisible group  of dimension $m$ over $\ol$.
 Suppose that
$V\times_{\spe K} \spe L\cong Y\times_{\spe \ol} \spe L$. We
denote by $G$ the local part of $Y$, by $J$ the corresponding
subgroup of $V$.

Then the following conditions are equivalent:

I There exists a $p$-divisible group $Z$ over $\ok$ such that
$J\cong Z\times_{\spe \ok} \spe K$.

II For a (finite flat commutative) group scheme $H/\ok$ we have
 $TH_\ol\approx (\ol/p^{l'}\ol)^m$; besides there exists a monomorphism
 $g:H_K\to\ke[p^{l'}]_{V,K}$.

III We have $TH_\ol\supset (\ol/p^{l'}\ol)^m$ (i.e. there exists
an embedding); there exists a monomorphism
 $g:H_K\to\ke[p^{l'}]_{V,K}$.

\end{pr}
\begin{proof}

 I $\implies$ II: We can take $H=[\ke
p^{l'}]_{Z,\ok}$. Then, according to part II of Theorem
\ref{dimgs}, $TH_\ol\approx (\ol/p^{l'}\ol)^m$.

II $\implies$ III: obvious.

Now we prove III $\implies$ I. We can assume that $H$ is local.

According to Theorem \ref{maingf}, there exists a morphism
$h:H_\ol\to \ke [p^t]_Y$ whose generic fibre is equal to $p^sg$.
We denote by $S$ the kernel of $h$ (i.e. $S\nde H$; $S_L=\ke h_L$;
see Proposition \ref{ray}); $T=H/S$. Since $h$ is defined over
$\ok$, so are $S$ and $T$.

We have an exact sequence $0\to S_\ol\to H_\ol\to T_\ol\to 0$;
besides $T_\ol\subset \ke [p^t]_G$. Since $p^sS=0$, we obtain $p^s
C(S)=0$; hence $p^s TS=0$. Since exact sequences of group schemes
correspond to exact sequences of tangent spaces (see part 3 of
Proposition \ref{dts}), we obtain that $TT\supset (\ol/p^t\ol)^m$.

Hence part 4 of Proposition \ref{dts} implies that $\ke
[p^t]_G=T$.

Since $T$ and the embedding are defined over $\ok$, according to
Proposition \ref{oldde}  we obtain that $G$ is defined over $\ok$
also.
\end{proof}

We also  prove the following (slight) generalization of Theorem
7.5.1 of \cite{02}.

\begin{theo}
Let $V$ be a $p$-divisible group  over $K$, let $Y$ be a
$p$-divisible group   over $\ol$.
 Suppose that
$V\times_{\spe K} \spe L\cong Y\times_{\spe \ol} \spe L$.

Then the following conditions are equivalent:

I. There exists a  $p$-divisible group $Z$ over $\ok$ such that
$V\cong Z\times_{\spe \ok} \spe K$.

II. For a (finite flat commutative) group scheme $H/\ok$ we have
 $TH_\ok\approx (\ok/p^{l'}\ok)^m$; besides there exists a isomorphism
 $g:H_K\to\ke[p^{l'}]_{V,K}$.
\end{theo}
\begin{proof}
 I $\implies$ II Again we can take $H=[\ke
p^{l'}]_{Z,\ok}$.

II $\implies$ I In \cite{02} the same statement was proved for $H$
equal to $\ke[p^{l'}]_U$ for some $p$-divisible $U$ over $\ok$. It
was noted also that it is sufficient to check the natural map
$i:H/\ke[p^s]_H\to H$ is a closed embedding.

We verify that $H$ satisfies this condition. Let $T$ denote the
local part of $H$. By part II of Theorem \ref{dimgs} $T$ is equal
to $\ke [p^{l'}]_F$ for some formal group $F/\ok$. Hence the
restriction of $i$ to $T$ is a closed embedding. Since the
exponent of $H$ equals $p^{l'}$, we have $i\ob(T)=T$. Since $H/T$
is etale, $i$ induces a closed embedding $H/T\to H/T$. Therefore
$i$ is a closed embedding itself.
\end{proof}

\subsection{Duality and semistable reduction of Abelian varieties}
\label{pdual}

 In order to reduce the descent question for Abelian
varieties to the study of certain formal groups we will need the
following statement.

Let $V$ denote an Abelian variety that has semistable reduction
over $L$, let $V'$ denote the dual variety. We denote by $V_p$ the
$p$-torsion of $V$ (as a $p$-divisible group over $L$), denote by
$V_f$ by the formal part of the $p$-torsion of $V$ (i.e. the part
corresponding to the formal group $F$ of the N\'eron model of $V$
over $\ol$). Obviously $V_f$ corresponds to the maximal
$p$-divisible subgroup $V_p$ that is a generic fibre of a local
scheme over  $\ol$. Hence the formal part is functorial with
respect to isogenies of Abelian varieties (one could also deduce
this fact from part 1 of the theorem below).

We note also that the formal part of an Abelian variety is
Galois-stable over its field of definition. Hence if $V$ is
defined over $K$ then $V_f$ is equal to $V_{fK}\times_{\spe K}\spe
L$ for a certain uniquely defined $p$-divisible subgroup
$V_{fK}\subset V_p$ over $K$.

\begin{theo}\label{dual}

 Let $V_m\subset V_f$ denote the multiplicative type
part of $V_p$ i.e. the part corresponding to the maximal
multiplicative type  subgroup of $F$. Let $V'_m$ denote the
multiplicative type part of $V'_p$.

1. The Weil pairing for $V_p$ and $V_p'$ induces the Cartier
duality of $V_p/V_f$ with $V'_m$.

2. Suppose also that $V$ is defined over $K$. Then $V$ has
semistable reduction over $K$ if and only if there exists a formal
$p$-divisible group $G/\ok$ such that the generic fibre of $G$ (as
a $p$-divisible group) is isomorphic to $V_{fK}$.
\end{theo}
\begin{proof}

1. We denote by $V_{fin}$, $V'_{fin}$ the finite parts of $V_p$
and $V_p'$ respectively, by $V_t$, $V'_t$ the toric parts of $V_p$
and $V_p'$ respectively (see the definitions in  \cite{gro}).
Obviously, $V_t'\subset V'_m$.

By Theorem 5.2 of \cite{gro}) $V'_{fin}$ is the Cartier dual of
 $V_p/V_t$, $V'_{t}$ is
the Cartier dual of
 $V_p/V_{fin}$. Hence for $D$ being the Cartier dual of $V_p/V_f$
 we have $V_t'\subset D\subset V'_{fin}$. Since $V_{fin}/V_f$ is
 the maximal possible \'etale factor of $V_{fin}$ (over $\ol$), $D/V_t'$ is the
 maximal possible multiplicative type subgroup of $V'_{fin}/V'_t$.
 Hence $D=V'_m$.

2. We define $V'_{fK}$ similarly to $V_{fK}$. Let $\al:V\to V'$ be
some polarization over $K$. Since the formal part is functorial,
$\al$ induces an isogeny $\beta:V_{fK}\to V'_{fK}$. Then
$\ke\beta$ is a finite $K$-group scheme; hence it corresponds to
some $S\nde G$ (see Proposition \ref{ray}). Therefore
$V'_{fK}\cong Z'\times_{\spe \ok} \spe K$, where $Z'=Z/S$ is a
formal $p$-divisible group over $\ok$.

We obtain that $V'_{fK}$  can be defined over $\ok$. Since the
multiplicative type part of an Abelian variety  is Galois-stable,
$V'_{mK}$ (defined similarly to $V'_{fK}$) is defined over $\ok$
also. Using the assertion of part 1,  we obtain that $V_{fK}$ and
$V_{pK}/V_{fK}$ are defined over $\ok$. Hence $V_{pK}$ is
Barsotti-Tate of echelon $2$ over $\ok$. By Proposition 5.13c  of
\cite{gro}, $V$ has semistable reduction over $K$.
\end{proof}

\begin{rema}
In \cite{gro} $V_{fin}$ was called the fixed part of $V_p$, $V_m$
was called the toroidal part of $V_p$, $V_{K,fin}$ was called the
effectively fixed part.
\end{rema}

\subsection{Criterion for good reduction of a potentially good reduction Abelian variety}

\begin{theo}\label{negr}
Let $V$ be be an Abelian variety  of dimension $m$ over $K$
 that has good reduction over $L$.
 Then $V$ has good reduction over $K$ if and only if
  for a certain (finite flat commutative) group scheme $H/\ok$ we have
   $TH_\ol\supset  (\ol/p^{l'}\ol)^m$  (i.e. there exists an embedding) and
  there exists a monomorphism
 $g:H_K\to\ke[p^{l'}]_{V,K}$.
\end{theo}
\begin{proof}
If $V$ has good reduction over $K$ then we can take $H=[\ke
p^{l'}]_{Y,\ok}$, where $Y_\ok$ is the N\'eron model of $V$ over
$\ok$.

Now we prove the converse implication. By  Proposition \ref{degr}
the formal part of $V_p$ is defined over $\ok$. By  part 2 of
Theorem \ref{dual} we conclude that $V$ has semistable reduction
over $K$. Hence $V$ has good reduction over $K$.
\end{proof}

In \cite{02} we proved a  version of the statement  where any
group scheme $H/\ok$ was allowed. We required the condition
$H_K\approx \ke[p^{l}]_{V,K}$.

\section{Finite criteria for semistable and ordinary reduction}

\subsection{Proof of Theorem A}

We prove Theorem A of the introduction. We adopt the notation of
Theorem \ref{dual}.

If $V$ has semistable reduction over $K$ then we can take $H=[\ke
p^l]_{G}$.

Now we prove the converse implication in Theorem A. We can assume
that $H$ is local.

 Since the finite part of an Abelian variety is Galois-stable
 over its field of definition,
we obtain that $V_{fin}$  is equal to $V_{finK}\times_{\spe K}\spe
L$ for a certain $p$-divisible group $V_{finK}$ over $K$.
According to Proposition 5.6 of \cite {gro}, the factor
$V_p/V_{fin}$ is the generic fibre of some \'etale $p$-divisible
group over $\ol$. Let $H_0$ denote the closed group scheme of $H$
that corresponds to the preimage of $V_{fin}$ in $H_L$ (see
Proposition \ref{ray}). We obtain that $(H/H_0)_L$ is isomorphic
to a generic fibre of a certain $\ol$-group scheme. Since $H$ is
local, Theorem \ref{maingf}  implies that $p^s(H/H_0)=0$ (in fact,
this statement is very easy to prove).

We denote by $H_1$ the maximal closed subgroup scheme of $H$ that
is killed by $p^{l'}$. The same tangent space argument as in the
proof of Proposition \ref{degr} shows that $TH_1\supset
(\ol/p^{l'}\ol)^m$.

We have ${H_1}\nde {H_0}$, whence ${H_0}_K\nde V_{fin K}$.
Therefore $TH_0\supset (\ol/p^{l'}\ol)^m$. We use the same
argument as in the proof of Theorem \ref{negr}. By Proposition
\ref{degr} the formal part of $V_p$ is defined over $\ok$. By part
2 of Theorem \ref{dual} we conclude that $V$ has semistable
reduction over $K$.

\subsection{Ordinary reduction: reminder}\label{ordred}

We recall that an Abelian variety (over $\ok$ or $\ol$) is called
an ordinary reduction one (or just ordinary) if the formal group
of its N\'eron module is finite height and of multiplicative type.

In particular, an ordinary variety has semistable reduction.

For example, a semistable reduction elliptic curve is either
ordinary or supersingular.

One can easily seen that the definition given is equivalent to the
usual ones. In \cite{mt} an Abelian variety was called ordinary
over $\ol$ if the connected component of $0$ of the reduction of
the N\'eron model of $A$ over $\ol$ is an extension of a torus by
an ordinary Abelian variety (over $\ovl$).

A good reduction Abelian variety that is ordinary is called a good
ordinary reduction one.

\subsection{Proof of Theorem B of the introduction}

\begin{proof}
I (1) $\implies $ (2)

We can take  $H=\ke [p^{l'}]_{F,\ok}$
 where $F$ is the formal group of the N\'eron model of $V$ (over $\ok$).

(2) $\implies $ (1)

According to Theorem \ref{negr} if such an $H$ exists then $V$ has
good reduction over $K$. The same argument as in the proof of
Theorem \ref{degr} shows that $\ke [p^t]_{F,\ol}$ is
multiplicative. Then $F$ is also multiplicative. In order to see
that one may pass to the duals and notice that a $p$-divisible
group $U$ is \'etale iff $\ke [p^t]_U$ is.

(2) $\implies $ (3)

Let $H'$ denote the Cartier dual of $H$. Since $H'$ is \'etale,
$H'$ is constant over some $A$  unramified over $\ok$. We take $M$
being the fraction field of $A$.

Since $H'_A$ is constant, $H_A$ is isomorphic to a
 $\sum_{1\le i\le r}\mu_{p^{n_i}}$ for some $n_i,r>0$.
We have $\mu_{p^{n_i},A}\cong \spe B$, where
$B=A[x]/((x+1)^{p^{n_i}}-1)$ and $J(B)=(x)$. Hence
$T\mu_{p^{n_i}}\approx A/p^{n_i}A$. Since $TH_{A\ol}\approx
(A\ol/p^{l'}A\ol)^m $, using the additivity of the tangent space
functor we obtain that $r=m$ and $n_i=l'$ for $1\le i\le r$. Thus
$H_M\cong (\mu_{p^{l'},M})^m$.

(3) $\implies $ (1)

We denote the ring of integers of $M$ by $A$. We  take $H_A=
(\mu_{p^{l'},A})^m$. Since $TH\approx (A/p^{l'}A)^m$, $H_A$
satisfies the conditions of (2) over $M$. Hence $V$ has good
reduction over $M$. Since $M/K$ is unramified,  $V$ has good
reduction over $K$ also.

II (1) $\implies $ (2)

If $V$ has  ordinary reduction, we can take  $H=\ke
[p^{l}]_{F,\ok}$.

 (2) $\implies $ (1)

 We apply Theorem A.
As in the proof of part I, we obtain that $\ke [p^t]_F$
  is multiplicative, hence $F$ is multiplicative.
We also note that (1) and (2) are equivalent to (2) with the
condition $(\ol/p^l\ol)^m\subset TH$  replaced by $TH_\ol\approx
(\ol/p^{l}\ol)^m$.

(2)$\iff$ (3)
 The proof (for $TH_\ol\approx (\ol/p^{l}\ol)^m$)
is the same as in part I.

\end{proof}

\begin{rema} 1. Passing to the duals, one may replace
a multiplicative  subgroup scheme of $\ke[p^l]_V$ by an \'etale
factorscheme in part II(3) of Theorem B.

2. (1)$\iff$(3) of Theorem B is quite similar to the finite
$l$-adic semistable reduction criteria of \cite{sz1} and
\cite{sz2}.
\end{rema}

\end{document}